 %
 %
 %
 %
 %
 %
 %
 %
\input amstex
\documentstyle{amsppt}

\hsize=4.75in
\vsize=8in

\hoffset=19mm
\font\small=cmr8

\def\sgap{\par\smallskip}
\def\newline{\hfil\break\noindent}

\def\cstar{C*}                \def\fidi{finite-dimen\-sional}
\def\rep{{repre\-sen\-tation}}   
\def\qft{quantum field theor}    \def\mo{\-mor\-phism}
\def\sps{{super\-se\-lec\-tion}} \def\hsp{Hilbert space}
\def\vn{von~Neumann}             \def\cex{conditional expectation}
\def\ca{Cuntz algebra}           \def\gca{generalized \ca}        

  \def\bb{{\Cal B}}  \def\hh{{\Cal H}}  \def\oo{{\Cal O}}
\def\CC{{\Bbb C}} \def\NN{{\Bbb N}}  \def\ZZ{{\Bbb Z}}
\def\eins{{\mathchoice {\hbox{\rm 1}\mskip-4mu\hbox{\rm l}} 
{\hbox{\rm 1}\mskip-4mu\hbox{\rm l}} {\hbox{\rm 1}\mskip-4.5mu\hbox{\rm l}} 
{\hbox{\rm 1}\mskip-5mu\hbox{\rm l}}}}
\def\sig{{\sigma}}    \def\a{\alpha}  \def\rho{\varrho}   
\def\id{\hbox{id}}   \def\comp{{\scriptstyle \circ}}    
\def\frac#1#2{{#1\over#2}}     \def\arr#1{\colon#1 \to}     
\def\red{_{\hbox{\small red}}}    \def\inv{^{-1}}     
\def\iro{_{\iota,\rho}}  \def\ro{_\rho}   \def\ioro{\iota\comp\rho} 
\def\text#1{\qquad\hbox{#1}\qquad}        \def\QED{\hfill $\qed$}

\rightheadtext {Generalized Cuntz algebras}
\leftheadtext {K.-H. Rehren}
\topmatter
\title
Generalized Cuntz algebras associated with subfactors
\endtitle
\author
Karl-Henning Rehren
\endauthor
\affil
II. Institut f\"ur Theoretische Physik, Univ.\ Hamburg (Germany)
\endaffil
\abstract
Various generalizations of Cuntz algebras and their relations to 
symmetry and duality are reviewed. New generalized Cuntz algebras are 
associated with a subfactor. A characteristic \hsp\ of basic invariants
(with respect to the generalized symmetry) within these algebras is 
discussed.
\endabstract
\endtopmatter
\document
\def\CDPR{1} \def\C{2} \def\DRdual{3} \def\DRwhy{4} \def\Iz{5} \def\Lo{6} 
\def\LH{7} \def\LR{8} \def\Pi{9}
\def\cite#1{\hbox{[#1]}}
\heading
1. Motivation
\endheading 
We consider an irreducible inclusion $A \subset B$ of properly infinite \vn\ 
factors, and denote by $\iota \in Hom(A,B)$ the inclusion map $A \ni a 
\mapsto a \in B$.
We assume the index to be finite. 
Then there is a conjugate endo\mo\ $\bar\iota\arr BA$ such that 
$\iota\comp\bar\iota\in End(B)$ is the canonical endo\mo\ of $B$ into $A$, 
and $\bar\iota\comp\iota\in End(A)$ is its dual \cite\Lo. 
We are interested in the irreducible subsectors $\rho$ of 
$\bar\iota\comp\iota$. 
If $\rho \prec \bar\iota\comp\iota$, then equivalently $\iota \prec \ioro$, 
that is, there exist intertwiners $\psi$ in $B$ which satisfy 
$$ \psi \cdot \iota(a) = \iota\rho(a) \cdot \psi \qquad\quad (a \in A) \; .
\eqno(1.1) $$
The existence of nontrivial such intertwiners, in turn, characterizes the
subsectors of $\bar\iota\comp\iota$.

In a specific \qft etical context, $A$ and $B$ may be thought of as a pair
of local algebras of observables and of charged fields, respectively.
Under standard physical assumptions (``unbroken symmetry'', cf.\ \cite\LR)
the dual canonical endo\mo\ $\bar\iota\comp\iota$ of $A$ extends to a 
localized endo\mo\ of the \cstar-algebra of observables, and 
can be identified with the restriction of the vacuum \rep\ of the field
algebra to the observables \cite\LR. The intertwiners (1.1) are then
field operators interpolating between the vacuum \rep\ of the observables
and a charged \rep\ corresponding to $\rho$. We shall therefore call these
intertwiners ``charged operators'' also in general. The charged operators 
form a \hsp\ $H$ of isometries within $B$. 

If $A \subset B$ has depth 2, then $A$ are the fixpoints under the action of 
a Hopf algebra (e.g., a group) $G$ on $B$, and the \hsp\ of charged operators 
for a given subsector $\rho$ has support projection $\eins_B$ in $B$ 
and carries a unitary \rep\ of $G$ \cite\LH. This establishes a correspondence 
between subsectors of $\bar\iota\comp\iota$ and \rep s of the symmetry. The 
duality problem amounts to the reconstruction of $G$ and $B$ from the 
knowledge of the relevant endo\mo s $\rho \in End(A)$ only 
\cite{\DRdual,\DRwhy}. 

One way to do so (if $G$ is a finite group) is to construct, for suitable 
$\rho$ of dimension $d$, the \ca\ $\oo(H)\cong\oo_d$ generated by a 
$d$-dimensional 
\hsp\ $H$ of isometries of unit support \cite\C, and to extend $A$ by 
$\oo(H)$ by postulating the commutation relations (1.1) for $\psi \in H$. 
One thus obtains a ``crossed product by the action of the endo\mo''. 
In order to recover $B$, however, a certain subalgebra of the \ca\ has to be 
identified with the subalgebra $\oo\ro$ of $A$ generated by the intertwiners 
$T \in A$ between powers of $\rho$, that is 
$$ T \cdot \rho^n(a) = \rho^m(a) \cdot T \qquad\quad (a \in A) \; . 
\eqno(1.2) $$
The unitary transformations of $H$ induce auto\mo s of $\oo(H)$. The 
identification of the subgroup $G \subset U(H)$ (the symmetry group) which has 
$\iota(\oo\ro)$ as its fixpoints within $\oo(H)$, constitutes the crucial 
step in the reconstruction problem \cite\DRdual. 

We want to understand the role of the charged operators for the position of 
$\iota(A)$ in $B$ in the case of arbitrary depth, i.e., not given by the f
ixpoints under a Hopf algebra. 
One could again consider the \cstar-subalgebra of $B$
generated by the charged operators, and study its relation with the 
\cstar-algebra $\oo\ro \subset A$ generated by all elements satisfying (1.2) 
for some $n,m$. 

But the latter algebra, considered as a subalgebra of $B$, is in general 
not contained in the former. It is therefore desirable to identify a 
larger \cstar-algebra $\oo\iro$ containing the charged operators as well as 
the invariants $\oo\ro$. In order to keep this algebra minimal, we require 
that its image under the unique normal \cex\ $\mu\arr BA$ gives exactly 
$\oo\ro$: 
$$ \mu(\oo\iro) = \oo\ro \text{and} \iota(\oo\ro) \subset \oo\iro \; .
\eqno(1.3) $$
The general idea is to consider the map $\mu\arr{\oo\iro}\oo\ro$ as a ``model
action'' for the generalized symmetry underlying $\mu\arr BA$.

This is the setting of our present analysis. We shall assume that $\iota$ and 
$\rho$ have finite index, and that both $\rho^n$ and $\ioro^n$, as $n$ 
varies, generate only finitely many irreducible subsectors (``rationality''). 
On the other hand, we drop the assumption that $\iota$ is contained in $\ioro$ 
(or $\rho\prec\bar\iota\comp\iota$), that is, possibly there are no charged 
operators in the proper sense. 
By rationality, however, $\iota$ will be contained in $\ioro^n$ for some $n$, 
so there will always be ``higher'' charged operators for $\rho^n$. 

Besides, $A$ need not be a subfactor of $B$, but rather $\iota$ may be any 
irreducible homo\mo\ between two properly infinite factors. Identifying
$A$ with its image $\iota(A)$ in $B$, the former situation will always be 
recovered.

With these data we associate two 
\cstar-algebras, namely $\oo\ro \subset A$ which is generated by all operators 
in $A$ which intertwine powers of $\rho \in End(A)$ as before, and $\oo\iro 
\subset B$ which is generated by all operators in $B$ which intertwine powers 
of $\rho$ (for the image of $A$ in $B$ like in (1.1)).

In order to analyze these \gca s, we study in both algebras the linear 
subspace of intertwiners between $\id_A = \rho^0$ and $\rho^n$, and identify 
within this subspace a generating \hsp\ of isometries (``skeleton space''). 
Surprisingly, in some exceptional cases the skeleton space is \fidi\ and 
has unit support. In these cases, it actually generates the entire \gca. 
In general, the skeleton space is in\fidi, but its support always
converges to unity in a natural Hilbert norm.

These observations give rise to embeddings of ordinary \ca s (generated by
the skeleton spaces) into \gca s. While {\it a priori\/} such embeddings are 
nothing peculiar,  the ones due to the skeleton space are of some relevance 
to the motivating duality problem, since the skeleton space of $\oo\iro$ 
collects the algebraically independent ``higher'' charged operators.

The \ca s generated by the skeleton spaces are in general dense within the 
\gca s (in an appropriate topology). In the exceptional cases with skeleton 
spaces of finite dimension, there is equality, that is, $\oo\iro$ may in fact 
turn out to be the ordinary \ca\ generated by its skeleton space.

As a byproduct, we give a simple characterization of the exceptional cases 
in terms of fusion matrices. 

There is a certain overlap with work done by Izumi \cite\Iz. In particular,
the $A_4$-example described in Sect.\ 3 can be found also in that paper, 
although considered from a different angle of view. 
We are not aware of a relation to another recent generalization of \ca s 
due to Pimsner \cite\Pi.

\heading
2. Generalized Cuntz algebras
\endheading
The \ca s $\oo_d$ were introduced in \cite\C\ as the \cstar-algebras 
generated by a (possibly countably-infinite) number $d \geq 2$ of orthogonal 
isometries $v_i$, subject to the relation $\sum_i v_iv_i^* = \eins$
if $d < \infty$. They were shown to be simple \cstar-algebras, depending only 
on $d$ up to iso\mo s. We are here mainly interested in concrete 
realizations of \ca s on a \hsp\ $\hh$. 
\proclaim{2.1 Definition} Let $H \subset \bb(\hh)$ be a \hsp\ of isometries,
that is $v^*w$ is a multiple of $\eins_\hh$ for $v,w \in H$. Then $\oo(H) 
\subset \bb(\hh)$ is the \cstar-algebra generated by the elements of $H$:
$$ \oo(H) := \left[\bigcup_{n,m\geq 0}H^m(H^n)^*
\right]^{\hbox{\raise-.8ex\hbox{$\scriptstyle C^*$}}} \subset \bb(\hh) \; .
$$
\endproclaim
\noindent The abelian case $\dim H = 1$ is of no interest.
If $H$ is separable, or if $2\leq \dim H < \infty$ and $H$ has unit support, 
then $\oo(H)$ is isomorphic with the \ca\ $\oo_{\dim H}$, by identification 
of $v_i$ with an orthonormal basis of $H$. 
In contrast, if $H$ is \fidi\ but of support $<\eins$, then 
$\oo(H)$ is known as a Cuntz-Toeplitz algebra which has a nontrivial ideal 
generated by $\sum_i v_iv_i^* - \eins$.

It was observed in \cite\DRdual\ in a new approach to duality for compact 
groups, that $\oo_\infty$ is not the appropriate object for duality theory.
The authors rather conceptualized the definition by viewing the generating
\hsp\ $H$ as an object of the category of \hsp s. This leads to a different 
generalization to the in\fidi\ case. Their general definition is 
\proclaim{2.2 Definition} Let $X$ be an object of a strict \cstar\ tensor
category. The spaces of arrows between tensor $($= monoidal\/$)$ 
powers of $X$ are 
denoted by $(X^{\times n},X^{\times m})$. These spaces are embedded into 
$(X^{\times n+1},X^{\times m+1})$ by taking the right tensor $($=
monoidal\/$)$ product with $\eins_X$. Then
$$ \oo_X := \left[\bigcup_{n,m\geq 0}(X^{\times n},X^{\times m})
\right]^{\hbox{\raise-.8ex\hbox{$\scriptstyle C^*$}}} \; .
$$
\endproclaim
{\it Example 2.1:} The category of \hsp s is a strict \cstar\ tensor category
where the arrows $(H_1,H_2)$ are given by the homo\mo s $(H_1,H_2) = 
Hom(H_1,H_2)$. For a \fidi\ \hsp\ $H$ of unit support, $\oo_H$ is isomorphic 
with $\oo(H)$ by the natural identification of $(H^{\otimes n},H^{\otimes m}) 
\subset \oo_H$ with $H^m(H^n)^* \subset \oo(H)$:
$$ \oo_H \cong \oo(H) \cong \oo_{\dim H} \; .
\eqno(2.1) $$
In contrast, for $H$ in\fidi, $\oo(H)$ is naturally embedded as a proper 
subalgebra into $\oo_H$. The reason is, in a nutshell, that $(H,H) \subset 
\oo_H$ is isomorphic with $\bb(H)$, while its intersection $HH^*$ with 
$\oo(H)$ contains only the compact operators.
\sgap
{\it Example 2.2:} The category of unitary \rep s $D \arr G\bb(H)$ of a 
compact group (or Hopf algebra) $G$ is a strict \cstar\ tensor category where 
the arrows $(D_1,D_2)$ are the intertwiners $Hom_G(H_1,H_2)$. 
If $D$ is a \rep\ of $G$ on the \fidi\ \hsp\ $H$ of isometries, then 
$(D^{\otimes n},D^{\otimes m})$ is naturally identified with the invariants 
under the induced action of $G$ on $H^m(H^n)^* \subset \oo(H)$, hence
$$ \oo_D \cong \oo(H)^G \; .
\eqno(2.2) $$
\sgap
Of particular interest for us is the case when $X = \rho$ is a unital 
endo\mo\ of some properly infinite \vn\ factor $A$. The endo\mo s of a 
properly infinite factor form a strict \cstar\ tensor category, where the 
arrows are the intertwiners, 
$$(\rho_1,\rho_2) = \{T \in A: T\rho_1(a) = \rho_2(a)T \quad\forall a \in A\}
              \qquad (\rho_1,\rho_2 \in End(A)) \; . 
\eqno(2.3) $$ 
The following definition is just a special case of Definition 2.2.
We assume $\rho$ to be proper, that is, not an auto\mo. 
\proclaim{2.3 Definition} Let $\rho$ be a proper unital endo\mo\ of $A$. 
Then
$$ \oo\ro := \left[\bigcup_{n,m\geq 0}(\rho^n,\rho^m)
\right]^{\hbox{\raise-.8ex\hbox{$\scriptstyle C^*$}}} \subset A \; .
$$
\endproclaim
\noindent Since $\rho$ maps $(\rho^n,\rho^m)$ into $(\rho^{n+1},\rho^{m+1})$, 
its action on $A$ preserves the subalgebra $\oo\ro$.
\sgap
{\it Example 2.3:} If $\rho$ is inner, i.e., $\rho(a) = \sum_i v_iav_i^*$
with isometries $v_i \in A$, then 
$$ \oo\ro = \oo(H) \; , 
\eqno(2.4) $$
where $H$ is the \hsp\ spanned by $v_i$ with inner product $(v,w) = v^*w$. 
As an endo\mo\ of $\oo(H)$, $\rho$ is the canonical inner endo\mo\ 
associated with $H$.
\sgap
We now introduce a new \gca, denoted $\oo\iro$, which is associated with an
irreducible homo\mo\ $\iota\arr AB$ between 
two properly infinite \vn\ factors $A$ and $B$, and a unital endo\mo\ $\rho$ 
of $A$ as before. The intertwiners for homo\mo s between properly infinite 
factors are given by 
$$(\sig_1,\sig_2) = \{T \in B: T\sig_1(a) = \sig_2(a)T \quad\forall a \in A\}
      \qquad (\sig_1,\sig_2 \in Hom(A,B)) \; . 
\eqno(2.5) $$ 
\proclaim{2.4 Definition} Let $\iota$ be an irreducible unital homo\mo\ of 
$A$ into $B$, and $\rho$ a proper unital endo\mo\ of $A$. Then
$$ \oo\iro := \left[\bigcup_{n,m\geq 0}(\ioro^n,\ioro^m)
\right]^{\hbox{\raise-.8ex\hbox{$\scriptstyle C^*$}}} \subset B \; . 
$$
\endproclaim
If $\mu$ is a \cex\ from $B$ onto $\iota(A)$, then $\mu$ maps
$(\ioro^n,\ioro^m)$ onto $\iota(\rho^n,\rho^m)$, hence 
$$ \mu(\oo\iro) = \iota(\oo\ro) \; .
\eqno(2.6) $$

The special case $B = A$, $\iota = \id_A$ yields $\oo\iro = \oo\ro$, of 
course. 
Our original motivation concerns the case $A \subset B$ and $\iota$ the 
inclusion map $A \ni a \mapsto a \in B$, and $\rho$ an irreducible subsector 
of $\bar\iota\comp\iota$. Then equivalently $\iota \prec \ioro$, that is, 
the linear subspace of charged operators in $B$
$$ H\iro = (\iota,\ioro) = \{\psi \in B: \psi a = \rho(a) \psi \quad 
\forall a \in A\}
\eqno(2.7) $$
is a non-trivial \hsp\ of isometries contained in $\oo\iro$. 
\sgap
{\it Example 2.4:}
Let $A \subset B$ irreducible be the fixpoints under a Hopf algebra $G$, and 
$\rho \prec \bar\iota\comp\iota$. 
Then $H\iro$ has unit support in $B$ and is stable under the action of $G$ 
and therefore carries a unitary representation $D$ of $G$ \cite\LH. 
Furthermore, as an endo\mo\ of $A$, $\rho$ is implemented within $B$ by 
$H\iro \subset B$,
$$\rho(a) = \sum_i\psi_i\;a\;\psi_i^* \qquad\quad \forall a \in A \; . 
\eqno(2.8) $$
One has  
$$ \oo\iro = \oo(H\iro) \cong \oo_{H\iro} \cong \oo_d 
\eqno(2.9) $$
where $d = d(\rho) = dim(H\iro) = dim(D)$. 
The \gca\ $\oo\ro \subset A$ equals the fixpoints of $\oo\iro \subset B$ under 
the Hopf algebra action. 
In view of Example 2.2, $\oo\ro$ can be identified with $\oo_D$:
$$ \oo\ro = \oo\iro \cap A = \left[\oo\iro\right]^G \cong \oo_D \; .
\eqno(2.10) $$
$\oo\ro$ may even be isomorphic with $\oo\iro$ itself, namely when $\rho = 
\bar\iota\comp\iota$ and $D$ is the regular \rep\ \cite{\CDPR,\Iz,\LH}; but in 
general it cannot be expected to be of the form $\oo(H)$ for any \hsp\ $H$.
\sgap
In general, the support of $H\iro$ is a nontrivial projection $<\eins_B$, 
and $\oo\iro$ is not generated by the \hsp\ $H\iro$. Because the dimensions of 
endo\mo s are multiplicative and additive, and because $\dim(H\iro)$ equals 
the multiplicity of $\iota$ within $\ioro$, one always has $\dim(H\iro) \leq 
d(\rho)$. But our assumptions even admit the possibility $H\iro = \{0\}$.

The previous examples are rather special illustrations. We emphasize that 
they do not exhaust the classes of \gca s of type $\oo\ro$ and $\oo\iro$.

\heading
3. The skeleton space
\endheading
We consider $\oo\iro$ as defined in Section 2. (The special case $B = A$, 
$\iota = \id_A$ covers also $\oo\ro$.) We assume $\iota$ and $\rho$ to be 
of finite dimension, that is $\iota(A) \subset B$ and $\rho(A) \subset A$ 
are subfactors of finite index. 

For each $n \in \NN_0$, $\oo\iro$ contains the \fidi\ subspace $h_n = 
(\iota,\ioro^n)$. These are \hsp s of isometries since $h_n^*h_n \subset 
(\iota,\iota) = \CC$. Being subspaces of $B$, $h_n$ can be multiplied, and 
$h_nh_m \subset h_{n+m}$. Accordingly, let $h'_n$ denote the span of all 
subspaces $h_{n_1} \cdots h_{n_r}$ with $1 \leq n_i < n$ and $\sum_i n_i = n$ 
within $h_n$, and $k_n$ the orthogonal complement of $h'_n$ in $h_n$. 
Trivially, $h_0 = \CC$, $k_1 = h_1 = H\iro$. Then
\proclaim{3.1 Lemma}  $x^*y = 0$ for $x \in k_n$, $y \in k_m$, provided 
$m \neq n$ and $n,m \neq 0$.
\endproclaim
\demo{Proof} Without loss of generality we may assume $m>n$. Then $x^*y 
\in h_{m-n}$.
It follows that $x \cdot x^*y \in h_nh_{m-n} \subset h'_m$, and by definition
this space is orthogonal to $k_m$, hence $y^*(x\cdot x^*y) = (x^*y)^*(x^*y) = 
0$. \QED
\enddemo
\proclaim{3.2 Definition} The skeleton space associated with $\oo\iro$ is 
the \hsp
$$ K\iro := \bigoplus_n k_n  \subset \oo\iro \; .
$$
For $B = A$, $\iota = \id_A$, we call $K\ro \equiv K\iro$ the skeleton space 
associated with $\oo\ro \equiv \oo\iro$.
\endproclaim
\noindent The obvious inclusions $H\iro \subset K\iro \subset \oo(K\iro) 
\subset \oo\iro$ hold. 

We introduce the power series
$$ h(t) = \sum_{n\geq 0} \dim(h_n) t^n  \text{and}
k(t) = \sum_{n>0} \dim(k_n) t^n 
\eqno(3.1) $$
as generating functionals for the dimensions of $h_n$ and $k_n$.
By construction, $h_n$ are generated by the subspaces $k_{n_1} \cdots k_{n_r}$
for all partitions $\sum n_i = n$. These are mutually orthogonal by
Lemma 3.1. A little combinatorics therefore yields
\proclaim{3.3 Lemma} The generating functionals for the dimensions of $h_n$ 
and $k_n$ are related by
$$ h(t) = \frac 1{1-k(t)} \qquad \Longleftrightarrow \qquad 
k(t) = 1 - \frac 1{h(t)} \; . 
$$ 
\endproclaim
Next, let us assume that the family of irreducible subsectors $\sigma_i$ of
$\ioro^n$ as $n$ varies, is finite (``rationality''). This property holds, 
e.g., if $\iota(A) \subset B$ is a subfactor of finite depth and $\rho$ is 
contained in some power of $\bar\iota\comp\iota$, that is, $\rho$ is a vertex 
in the even subgraph of the principal graph associated with the subfactor. The 
assumption is also satisfied if $B = A$, $\iota=\id$, and $\rho$ is an object 
of a category of endo\mo s with only finitely many irreducible objects, e.g., 
the \sps\ sectors of a ``rational'' \qft y. We then have a finite fusion 
matrix $(N\iro)_{ij}$ given by the non-negative multiplicities of $\sigma_i$ 
within $\sigma_j\comp\rho$, where $\sig_i$ range over all irreducible 
subsectors of $\ioro^n$ and one of the indices ($i = \iota$) corresponds to 
$\sigma_\iota \equiv \iota$. Since $\dim(h_n) = [(N\iro)^n]_{\iota\iota}$, one 
obtains the formula for the generating functional
$$ h(t) = [(1-tN\iro)\inv]_{\iota\iota} \; , 
\eqno(3.2) $$
from which $k(t)$ can be obtained with the help of the Lemma 3.3. Explicitly,
if $P$ is the projection matrix onto the $\iota$-direction and $Q = \eins-P$,
one has
$$ \dim(k_n) = [N\iro Q N\iro \dots N\iro Q N\iro]_{\iota\iota} 
\equiv [N\iro (Q N\iro Q)^{n-2} N\iro]_{\iota\iota} \; .
\eqno(3.3) $$
Clearly, $\dim(K\iro) = k(1)$.

It is well known that the largest eigenvalue of $N\iro$ is given by the 
dimension $d(\rho)$ of $\rho$, and the eigenvector is the Frobenius 
eigenvector $F_i = d(\sig_i)$. It follows that the radius of convergence of 
the generating functional $h(t)$ is $1/d(\rho) < 1$.
\sgap
{\it Example 3.1:} Let $\rho \in End(A)$ be implemented by a \hsp\ $H \in B$. 
Then $h_1 = H$, $h_n = H^n = h'_n$ and $k_1 = H$, $k_n = \{0\}$ ($n>1$), hence 
$K\iro = H$. This yields $h(t) = (1-dt)\inv$ where $d = \dim(H) = d(\rho)$, 
and $k(t) = d\,t$. As mentioned before, $\oo\iro = \oo(K\iro) = \oo(H)$. 
Extending $\rho$ to an inner endo\mo\ $\sig$ of $B$, we have $K_\sig = K\iro$ 
and $\oo_\sig = \oo\iro = \oo(H)$.
\sgap
{\it Example 3.2:} On the other hand, the algebra $\oo\ro$ associated with 
$\rho 
\in End(A)$ from Example 3.1 will in general have a skeleton space of infinite 
dimension. If, for instance, $A$ are the fixpoints of $B$ under a non-abelian 
symmetry group $G$ and $\rho$ corresponds to a \rep\ $D$ of $G$, then this 
space is determined by the \rep\ theory of $G$; namely a basis of $k_n \subset 
K\ro$ corresponds to the independent $G$-invariant tensors within 
$D^{\otimes n}$. We display the generating functional $k(t)$ for $G = S_3$ 
and $\rho$ corresponding to the two-dimensional \rep: $k(t) = t^2/(1-t-t^2)$.
\sgap
{\it Example 3.3:} Let $\iota(A) \subset B$ be a subfactor with principal 
graph $A_4$. Then $\bar\iota\comp\iota \cong \id_A \oplus \rho$ with 
$d(\rho) = \frac{1+\sqrt5}2$. The relevant fusion rules are $\ioro \cong \iota 
\oplus \a$ where $\a$ is an iso\mo\ of $A$ with $B$, and $\a\comp\rho \cong 
\iota$. The above prescription yields
$$ h(t) = [(1-tN\iro)\inv]_{\iota\iota} = \frac 1{1-t-t^2} \text{and}
k(t) = t+t^2 \; . 
$$
In other words, $k_1$ and $k_2$ are one-dimensional, and all higher $k_n$ are
trivial. 
\sgap
{\it Example 3.4:} Consider now $\oo\ro$ with $\rho$ as in Example 3.3. In 
this case, the relevant fusion rules are $\id\comp\rho = \rho$ and $\rho^2 
\cong \id \oplus \rho$, yielding
$$ h(t) = [(1-tN\ro)\inv]_{00} = \frac{1-t}{1-t-t^2} \text{and}
k(t) = \frac{t^2}{1-t} \; . 
$$
This shows that $K\ro$ is in\fidi.

We can look at the Examples 3.3, 3.4 more explicitly. Namely, due to the 
``Lee-Yang'' fusion rules $\rho^2 \cong \id \oplus \rho$, there are isometries 
$v_1 \in (\rho,\rho^2)$ and $v_2 \in (\id,\rho^2)$ in $A$ satisfying $v_1v_1^*
+ v_2v_2^* = \eins_A$. Actually, this pair of isometries generates $\oo\ro$, 
so $\oo\ro$ is isomorphic with the ordinary \ca\ $\oo_2$. But only $v_2$ is 
contained in $K\ro$.

To study $\oo\iro$, one may choose $\a$ within its inner equivalence class 
such 
that $\a\comp\rho = \iota$. It follows that $\psi_i = \a(v_i)$ are isometries 
$\psi_1 \in (\iota,\ioro)$ and $\psi_2 \in (\a\comp\rho,\a\comp\rho^3) = 
(\iota,\ioro^2)$ {\it a forteriori}. This pair of isometries is a basis of 
$K\iro$. Trivially, $\iota(\oo\ro)$ is (properly) contained in $\oo\iro$, 
while 
we shall see below that $\oo\iro$ is in fact generated by $\psi_i$, that is 
$\oo\iro = \a(\oo\ro)$. Thus we have the proper inclusion
$$ \iota(\oo\ro) \subset \oo\iro = \a(\oo\ro) \; . 
\eqno(3.4) $$ 
We see in these examples that the skeleton space is not an invariant under 
iso\mo s of \gca s, but depends on the presentation in terms of $\iota$ and 
$\rho$.

\heading 
4. Density
\endheading 
We are concerned with the question how the skeleton space $K\iro$ and the 
\ca\ $\oo(K\iro)$ generated by this \hsp\ of isometries are embedded into
$\oo\iro$. 
Before we consider the general case, we return to the Examples 3.3, 3.4 
associated with the $A_4$ inclusion, discussed in the previous section. 
It exhibits the general mechanism in a most transparent manner.

We therefore assume, as in Example 3.3, $\iota \in Hom(A,B)$ with the fusion 
rules $\bar\iota\comp\iota \cong \id \oplus \rho$, $\iota\comp\rho \cong \iota 
\oplus \a$, $\a\comp\rho = \iota$. Let $X \in (\ioro^n,\ioro^m) \in \oo\iro$. 
Then $X$ is a linear combination of operators $T_aT_b^*$ and $S_cS_d^*$ with 
isometries $T_a \in (\iota,\ioro^m), T_b \in (\iota,\ioro^n)$ and $S_c \in 
(\a,\ioro^m), S_d \in (\a,\ioro^n)$. Since $(\iota,\ioro^n) = h_n$ and 
$(\a,\ioro^n) \subset (\a\comp\rho,\ioro^{n+1}) = h_{n+1}$, and since in turn 
$h_n$ are generated by $k_1$ and $k_2$ which span $K\iro$, we conclude that 
$X$ is in the \ca\ generated by $K\iro$, that is, $\oo\iro = \oo(K\iro) 
\cong \oo_2$.

The situation is different for the Example 3.4. The endo\mo\ $\rho$ is the same
as before, but $\oo\ro$ is generated by intertwiners between powers of $\rho$ 
only. So let $X \in (\rho^n,\rho^m) \in \oo\ro$. Similar as before, $X$ is a 
linear combination of operators $R_aR_b^*$ and $U_cU_d^*$ with isometries 
$R_a \in (\id,\rho^m) = h_m, R_b \in (\id,\rho^n) = h_n$ and $U_c \in 
(\rho,\rho^m), U_d \in (\rho,\rho^n)$. This time, the latter can not be 
directly identified with elements of $K\ro$. But $(\rho,\rho^n)$ is embedded 
into $(\rho^2,\rho^{n+1})$ and $\rho^2 \cong \id \oplus \rho$ is decomposed 
with the pair of isometries $v_1$ and $v_2$. Thus we may write
$$ U_cU_d^* = U_cv_2v_2^*U_d^* +  U_cv_1v_1^*U_d^* \; ,
\eqno(4.1) $$
where $U_cv_2 \in h_{m+1}, U_dv_2 \in h_{n+1}$. Hence the first term in the 
decomposition is in the \ca\ $\oo(K\ro)$. The second term is of the form 
$\tilde U_c \tilde U_d^*$ with $\tilde U_c \in (\rho,\rho^{m+1}), \tilde U_d 
\in (\rho,\rho^{n+1})$, and can be treated iteratively like $U_cU_d^*$ before. 
Thus we end up with a decomposition
$$ U_cU_d^* = \sum_{r=0}^R (U_cv_1^rv_2)(U_dv_1^rv_2)^* +
   (U_cv_1^{R+1})(U_dv_1^{R+1})^* 
\eqno(4.2) $$
where the terms in the sum are in $h_{m+r+1}h_{n+r+1}^*$, hence are generated
by $K\ro$. We shall show that the remainder term converges to zero in an 
appropriate topology (but not in the operator norm), as $R$ increases.

In order to do so, we recall some definitions pertaining to the general case.
For $\sig$ a unital homo\mo\ of $A$ into $B$ of finite dimension, the standard 
left-inverse $\phi_\sig$ is the positive unital map $\sig\inv\comp E_\sig$ 
where $E_\sig$ is the minimal conditional expectation of $B$ onto its 
subfactor $\sig(A)$. Explicitly, $\phi_\sig$ is of the form $\phi_\sig(x) = 
R^*\bar\sig(x)R$ where $\bar\sig$ is a conjugate endo\mo\ of $\sig$ and $R$ 
is an isometry in $(\id_A,\bar\sig\sig)$. By construction, $\phi_\sig$ 
satisfies $\phi_\sig(\sig(x)y\sig(z)) = x\phi_\sig(y)z$ and in particular 
$\phi_\sig\comp\sig = \id_A$. Standard left-inverses are related by the 
``intertwining property''
$$ d(\sig_1)\phi_{\sig_1}(xT) = d(\sig_2)\phi_{\sig_2}(Tx) \text{if}
T \in (\sig_1,\sig_2) \; .
\eqno(4.3) $$

We consider the state on $\oo\ro$ given by 
$$ \varphi\ro = \lim_{N\rightarrow\infty}\phi\ro^N \; . 
\eqno(4.4) $$ 
$\phi\ro$ maps $(\rho^n,\rho^m)$ into $(\rho^{n-1},\rho^{m-1})$ as long as 
$n,m > 0$. Hence, its repeated application takes $(\rho^n,\rho^n)$ eventually 
to the scalars $(\id_A,\id_A) = \CC$ where further application of $\phi\ro$ 
is trivial. On the other hand, if $n \neq m$, and without loss of generality 
$n < m$, then $(\rho^n,\rho^m)$ is mapped by $\phi\ro^n$ to the \hsp\
$h_{m-n}$. This space is mapped by $\phi\ro$ onto itself, and it can be shown 
that the restriction of $\phi\ro$ to $h_k$ ($k>0$) is a bounded map with norm 
$\leq d(\rho)\inv < 1$. Hence the limit $\varphi\ro = \lim_{N\rightarrow\infty}
\phi^N$ stabilizes on the homogeneous part of $\oo\ro$ and converges to zero 
on the inhomogeneous part of $\oo\ro$. In other words, $\varphi\ro$ is a 
state on $\oo\ro$ which is invariant under the auto\mo\ group on $\oo\ro$ 
defined by $\a_t(x) = e^{i(n-m)t}x$ if $x \in (\rho^n,\rho^m)$. The 
intertwining property (4.3) implies
$$ \varphi\ro(xy) = d(\rho)^{n-m} \varphi\ro(yx) \text{for} x \in 
(\rho^n,\rho^m) \; ,
\eqno(4.5) $$
that is, $\varphi\ro$ is actually a KMS state for this auto\mo\ group with 
temperature $2\pi/\ln d(\rho)$.

The intertwining property allows to compute the expectation values
$$ \varphi\ro(TT^*) = \frac{d(\tau)}{d(\rho)^n} 
\eqno(4.6) $$
whenever $T \in (\tau,\rho^n)$ is an isometry and $\tau\prec\rho^n$ is 
an irreducible subsector of $\rho^n$.

Similarly, we consider the state on $\oo\iro$ given by 
$$ \varphi\iro = \varphi\ro\comp\mu \; ,
\eqno(4.7) $$
where the unique \cex\ $\mu$ of $B$ onto $A$ coincides with the left-inverse 
$\phi_\iota$. Since $\mu$ maps $(\ioro^n,\ioro^m)$ onto $(\rho^n,\rho^m)$, it 
maps $\oo\iro$ onto $\oo\ro$. Similarly as before, one has
$$ \varphi\iro(VV^*) = \frac{d(\sig)}{d(\iota)d(\rho)^n} 
\eqno(4.8) $$
whenever $V \in (\sig,\ioro^n)$ is an isometry and $\sig\prec\ioro^n$ is 
an irreducible subsector of $\ioro^n$.

The faithful state $\varphi\iro$ gives rise, by the GNS construction, to the 
\hsp\ $\hh\iro$ which is the closure of the pre \hsp\ $\oo\iro$.
Specializing to $B = A$, $\iota = \id_A$ as before, the state $\varphi\ro$
gives rise to the GNS \hsp\ $\hh\ro \supset \oo\ro$.

Let us now return to our Example 3.4. We were decomposing an element of 
$\oo\ro$ 
into a part generated by the skeleton space $K\ro$, and a remainder term 
$(U_cv_1^{R+1})(U_dv_1^{R+1})^*$. This term converges to zero when considered 
as an element of the GNS \hsp\ $\hh\ro$, namely 
$$ \eqalign{\Vert(U_cv_1^{R+1})(U_dv_1^{R+1})^*\Vert_{\varphi\ro}^2 & = 
\varphi\ro[(U_dv_1^{R+1})(U_cv_1^{R+1})^*(U_cv_1^{R+1})(U_dv_1^{R+1})^*]
\cr & = \varphi\ro[(U_dv_1^{R+1})(U_dv_1^{R+1})^*] 
\cr & = d(\rho)^{-n-R-1} \; . \cr} 
\eqno(4.9) $$
by (4.6). It is, however, crucial that this is no convergence in the 
\cstar-algebra $\oo\ro$, since every remainder term is a partial isometry, 
and hence has unit norm. 

The example already illustrates the general scheme. In fact, we have
\proclaim{4.1 Lemma} Let $E_n = \sum_\a v_\a v_\a^*$ be the support 
projection of the \hsp\ $k_n$ where $\{v_\a\}$ is an orthonormal basis of 
$k_n \subset K\iro$. Then the expectation value $\varphi\iro(\sum_n E_n)$ 
converges to 1.
\endproclaim
\demo{Proof} 
Since $v_i \in (\iota,\ioro^n)$, one has $\varphi\iro(v_iv_i^*) = 
d(\rho)^{-n}$ due to (4.8), and $\varphi\iro(E_n) = \dim(k_n) d(\rho)^{-n}$. 
It follows that the sum $\sum_n \varphi\iro(E_n)$ converges to the value of 
the generating functional $k(t)$ at $t = 1/d(\rho)$ while to evaluate $k(t)$ 
at $t$ below $1/d(\rho)$ amounts to suppress the higher contributions to the
sum. But since $h(t)$ diverges as $t \nearrow 1/d(\rho)$, we conclude by Lemma 
3.3 that $k(1/d(\rho)) = 1$. \QED
\enddemo
\proclaim{4.2 Corollary} $\sum_n E_n$ converges to $\eins_B$ within the 
\hsp\ $\hh_{\iota,\rho}$. If $K\iro$ is \fidi, then it has unit support, 
$\sum_n E_n = \eins_B$.
\endproclaim 
\demo{Proof} 
Obvious, since $E_n$ are orthogonal projections by Lemma 3.1. \QED
\enddemo
\proclaim{4.3 Corollary} $\dim(H\iro) \leq d(\rho) \leq \dim(K\iro)$. 
Equality in either of these bounds implies equality in the other one,
and holds if and only if either of the following is true: \newline
{\rm (i)} $K\iro = k_1 \equiv H\iro$. \newline
{\rm (ii)} $H\iro$ has unit support. \newline
{\rm (iii)} $\rho$ is implemented by $K\iro$. \newline
{\rm (iv)} $\rho$ is implemented by $H\iro$.
\endproclaim
\demo
{Proof} The first bound holds because $\dim(H\iro)$ is the multiplicity 
of $\iota$ within $\ioro$, with equality if and only if (iv) is true. 
The second bound holds since every basis isometry $v_\a \in k_n$ contributes 
at most $1/d(\rho)$ to the sum $\sum_{n,\a}\varphi\iro(v_\a v_\a^*)$
$= 1$, with equality if and only if each contribution is exactly $1/d(\rho)$, 
which in turn is equivalent to (i). 

Equivalence of (i)--(iv) is seen as follows:
if $K\iro = H\iro$ then it is \fidi, hence it has unit support, hence 
$K\iro = H\iro$ implements $\rho$. Conversely, if $H\iro$ implements $\rho$, 
then its support must be unity, and if $H\iro$ has unit support, then $K\iro$ 
cannot be larger than $H\iro$. If $K\iro$ implements $\rho$, then $K\iro 
\subset (\iota,\ioro) = H\iro$. In either case, $K\iro = H\iro$ follows. \QED
\enddemo
In a similar way as Lemma 4.1, one obtains
\proclaim{4.4 Lemma} Every element of $(\ioro^n,\ioro^m)$ can be approximated
$($in the Hilbert norm of $\hh_{\iota,\rho})$ by sums of operators in 
$h_{m+r} h_{n+r}^*$.
\endproclaim
\demo{Proof} 
Every $x \in (\ioro^n,\ioro^m)$ is a sum of operators $T_aT_b^*$ 
and $S_cS_d^*$ with isometries $T_a \in (\iota,\ioro^m), T_b \in 
(\iota,\ioro^n)$ and $S_c \in (\sig_j,\ioro^m), S_d \in (\sig_j,\ioro^n)$
where $\sig_j$ ranges over all irreducible subsectors contained in $\ioro^n$, 
$n \in \NN$, which are different from $\iota$. The terms $T_aT_b^*$ are of the 
desired form. The latter are treated as follows. $(\sig_j,\ioro^n)$ is 
a subspace of $(\sig_j\comp\rho,\ioro^{n+1})$. According to the fusion rules 
$\sig_j\comp\rho \cong (N\iro)_{\iota j}\iota \oplus \sum_{i\neq \iota} 
(N\iro)_{ij} \sig_i$ we may rewrite 
$$ S_cS_d^* = \sum_\a S_c w_\a w_\a^* S_d^* + 
\sum_{i,\a} S_c v_{i,\a}v_{i,\a}^*S_d^* 
\eqno(4.10) $$
where $\{w_\a\}$ is an orthonormal basis of $(\iota,\sig_j\comp\rho)$ and 
$\{v_{i,\a}\}$ are orthonormal bases of $(\sig_i,\sig_j\comp\rho)$. The former 
terms are in $h_{m+1}h_{n+1}^*$, while the latter are of the form $\tilde S_c
\tilde S_d^*$ with isometries $\tilde S_c \in (\sig_i,\ioro^{m+1}), \tilde S_d 
\in (\sig_i,\ioro^{n+1})$ ($i \neq \iota$). Repeating this step $R$ times, we 
obtain a sum of terms in $h_{m+r}h_{n+r}^*$ ($r\leq R$) and a remainder 
$r_R$ of the form
$$ r_R = S_c\left(\sum v_1 \dots v_R v_R^* \dots v_1^*\right)S_d^* 
\eqno(4.11) $$
where $v_r = v_{i_r,\a_r}$ and the sum extends over all $i_r \neq \iota$ and 
over orthonormal bases $\{v_{i_r,\a_r}\}$ of 
$(\sig_{i_r},\sig_{i_{r-1}}\comp\rho)$. The Hilbert norm of the remainder is, 
by (4.8),
$$ \Vert r_R \Vert^2_{\varphi\iro} = \varphi\iro(r_R^*r_R)
= \sum \varphi\iro\left[S_d v_1 \dots v_R v_R^* 
\dots v_1^*S_d^*\right] = \sum\frac{d_{i_R}}{d(\rho)^{n+R}} \; .
\eqno(4.12) $$
The number of terms in this sum with fixed sector $i_R = i$ equals 
$(N\red^R)_{ij}$ where $N\red = Q N\iro Q$ is the fusion matrix 
with the line and column corresponding to $\iota$ deleted. 
This yields the estimate
$$ \Vert r_R \Vert_{\varphi\iro}^2 \leq const. \cdot
\left(\frac{\Vert N\red \Vert}{d(\rho)}\right)^R \; .
\eqno(4.13) $$
The norm of the reduced fusion matrix is strictly smaller than the norm 
$\Vert N\iro \Vert$ $= d(\rho)$ of the full fusion matrix, because the 
Frobenius eigenvector has a nonzero component in the $\iota$-direction. 
It follows that (4.13) converges to zero as $R$ increases. \QED
\enddemo
\proclaim{4.5 Corollary} The following statements are equivalent. \newline
{\rm (i)} $\oo\iro$ is generated by its skeleton space $K\iro$, that is 
$\oo\iro = \oo(K\iro)$. \newline
{\rm (ii)} $K\iro$ is \fidi. \newline
{\rm (iii)} The ``reduced fusion matrix'' $N\red = Q N\iro Q$ $($cf.\
Eq.\ $(3.3))$ is nilpotent.
\endproclaim
\demo{Proof} 
For the equivalence (i) $\Leftrightarrow$ (iii), consider the 
expansion discussed in the proof to Lemma 4.4. Since every term in the
remainder $r_R$ is a partial isometry of norm 1, the expansion converges 
within $\oo\iro$, that is in norm, if and only if the remainder vanishes for 
sufficiently large $R$. This is equivalent to nilpotency of $N\red$. 
The equivalence (ii) $\Leftrightarrow$ (iii) is obvious from Eq.\ (3.3). \QED
\enddemo

\heading
5. Discussion
\endheading
We have introduced a new class of \gca s $\oo\iro$ motivated by the analysis
of charged operators, that is, operators in an ambient algebra $B$ which 
intertwine a given endo\mo\ $\rho$ of a subalgebra $A$ with the identity
$\id_A$. 

We have studied the embedding of ``higher'' charged operators, that is,
intertwiners in $B$ between the identity and powers of $\rho \in End(A)$, into
the associated \gca. An interesting invariant object turns out to be the 
``skeleton space'' $K\iro$ within $\oo\iro$. In the fixpoint example, the 
skeleton space of $\oo\iro$ is just the \hsp\ $H\iro$ of charged operators 
which implements $\rho$ and carries a \rep\ of the Hopf algebra $G$, while 
the skeleton space of $\oo\ro$ collects a system of algebraically independent 
$G$-invariants within tensor powers of the \rep\ on $H\iro$. In this 
sense, we are setting up a generalized theory of invariants.

The skeleton space gives rise to an embedding $\oo(K\iro) \subset \oo\iro$. 
The embedding of ordinary \ca s into generalized ones itself is not so much 
surprising in view of the elementary fact that ordinary \ca s generated by
$d$ isometries can be embedded into \ca s generated by less than $d$ isometries
(e.g., if $\oo_2$ is generated by $v_1$ and $v_2$, then $v_1^rv_2$ ($r = 0, 
\dots d-2$) and $v_1^{d-1}$ span a \hsp\ of dimension $d$ with unit support
within $\oo_2$). But the simple criterion for equality of $\oo\iro$ with 
$\oo(K\iro)$ in terms of the fusion rules (Corollary 4.5) seems 
interesting, possibly leading to a classification of these exceptional cases.
The (rather weak) convergence of the support of the \hsp\ $K\iro$ in the 
general case may be irrelevant for the properly \cstar\ aspects, but in view 
of the motivating problem in which the states $\varphi\iro$ and $\varphi\ro$
are natural tools, we expect it to be of interest for the study of
generalized symmetry.

The criterion whether $\oo\iro$ is an ordinary \ca\ distinguishes only the 
case when the generating \hsp\ coincides with the skeleton space $K\iro$. 
The Example 3.4 shows that possibly $\oo\iro = \oo(H)$ with some 
$H \neq K\iro$. In that example, $K\iro$ is in\fidi\ while $\dim(H) = 2$. 
The skeleton space therefore does not provide a characterization of the 
full iso\mo\ class of the new \gca s.

However, the skeleton spaces and the generating functionals $k(t)$, as well 
as the relative position of $\iota(\oo\ro)$ in $\oo\iro$ are invariants under 
inner conjugations and outer transformations $(A,B,\iota,\rho) \mapsto 
(\hat A,\hat B,\beta\comp\iota\comp\alpha\inv,\alpha\comp\rho\comp\alpha\inv)$
with $\alpha\in Iso(A,\hat A)$ and $\beta\in Iso(B,\hat B)$.
These data, as $\rho$ ranges over the even vertices of the principal graph
determined by $\iota$, therefore constitute invariant information about the 
underlying paragroup.

\heading
References 
\endheading

\def\CMP#1{Com\-mun.~Math.~Phys.\ {\bf #1}}

\def\IM#1{Invent.~Math.\ {\bf #1}}
\def\JFA#1{J.~Funct.~Anal.\ {\bf #1}}

\def\RMP#1{Rev.~Math.~Phys.\ {\bf #1}}
\item{[1]} T.~Ceccherini, S.~Doplicher, C.~Pinzari, J.E.~Roberts: {\it A 
generalization of the \ca s and model actions}, \JFA{125}, 416--437 (1994).
\item{[2]} J.~Cuntz: {\it Simple \cstar-algebras generated by isometries},
\CMP{57}, 173--185 (1977).
\item{[3]} S.~Doplicher, J.E.~Roberts: {\it A new duality theory for
compact groups}, \IM{98}, 157--218 (1989).
\item{[4]} S.~Doplicher, J.E.~Roberts: {\it Why there is a field
algebra with a compact gauge group describing the \sps\ structure in
particle physics}, \CMP{131}, 51--107 (1990).
\item{[5]} M.~Izumi: {\it Subalgebras of infinite \cstar-algebras with
finite Watatani indices. I. \ca s}, \CMP{155}, 157--182 (1993).
\item{[6]} R.~Longo: {\it Index of subfactors and statistics of quantum 
fields.\ I}, Commun.~Math.~Phys.\ {\bf 126}, 217--247 (1989).
\item{[7]} R.~Longo: {\it A duality for Hopf algebras and for
subfactors.\ I}, Commun.~Math.\ Phys.\ {\bf 159}, 133--150 (1994).
\item{[8]} R.~Longo, K.-H.~Rehren: {\it Nets of subfactors}, \RMP{7},
567--597 (1995).
\item{[9]} M.~Pimsner: {\it A class of \cstar-algebras generalizing both
Cuntz-Krieger algebras and crossed products by $\ZZ$}, Univ.\ of Pennsylvania
preprint (1996).

\enddocument
\bye